\documentclass[11pt]{article}
\usepackage{flafter,amsmath,amssymb,latexsym,psfrag,graphicx,color,indentfirst}
\usepackage[margin=2.5cm]{geometry}

\usepackage[colorlinks,
linktocpage=true,
linkcolor=black,
citecolor=black,
urlcolor=black,
anchorcolor=black
]{hyperref}
\usepackage{setspace}
\usepackage[numbers,sort&compress]{natbib}
\usepackage{appendix}
\usepackage{multirow}
\usepackage{caption}
\usepackage{subfigure}
\usepackage{graphicx}
\usepackage{tabularx}
\usepackage{multicol}
\usepackage{multirow}
\usepackage{booktabs}
\usepackage{comment}
\usepackage{epsfig}
\usepackage{epstopdf}
\usepackage{amsthm,amsmath,amssymb}
\usepackage{mathrsfs}
\usepackage{ulem}
\graphicspath{{figures/}}

\newtheorem{theorem}{Theorem}[section]
\newtheorem{defn}[theorem]{Definition}
\newtheorem{lemma}[theorem]{Lemma}
\newtheorem{remark}[theorem]{Remark}

\newtheorem{prop}[theorem]{Proposition}

\makeatletter

\@addtoreset{figure}{section}
\makeatother

\begin{document}
\setlength\arraycolsep{2pt}
\date{\today}

\title{On the convergence of the no-response test for the heat equation}
\author{Shiwei Sun$^{1}$,\, Gen Nakamura$^{2,\,3}$,\, Haibing Wang$^{1,}$\footnote{Corresponding author, E-mail: hbwang@seu.edu.cn}\\
\\$^1$School of Mathematics, Southeast University, Nanjing 210096, P.R. China\;\;
\\$^2$Department of Mathematics, Graduate School of Science,
\\Hokkaido University, Sapporo 060-0810, Japan
\\$^3$Research Center of Mathematics for Social Creativity,
\\Research Institute for Electronic Science,
\\Hokkaido University, Sapporo 060-0812, Japan
}

\maketitle
\begin{abstract}
Domain sampling methods called the range test (RT) and no-response test (NRT), and their duality are known for several inverse scattering problems and an inverse boundary value problem for the Laplace operator (see Section 1 for more details). In our previous work \cite{Sun-2024-40}, we established the duality between the NRT and RT, and demonstrated the convergence of the RT for the heat equation. We also provided numerical studies for both methods. However, we did not address the convergence for the NRT. As a continuation of this work, we prove the convergence of the NRT without using the duality. Specifically, assuming there exists a cavity $D$ inside a heat conductor $\Omega$, we define an indicator function $I_{NRT}(G)$ for a prescribed test domain $G$, where $\overline G\subset\Omega$ (i.e., $G\Subset\Omega$). By using the analytical extension property of solutions to the heat equation with respect to the spatial variables, we prove the convergence result given as $I_{NRT}(G)<\infty$ if and only if $\overline{D}\subset \overline{G}$, provided that the solution to the heat equation cannot be analytically extended across the boundary of the cavity. Thus, we complete the theoretical study of both methods. Here the analytic extension of solutions does not require the property that the solutions are real analytic with respect to the space variables. However, for the proof of the mentioned convergence result, we fully use this property.



\end{abstract}

\section{Introduction}\label{Introduction}
\setcounter{equation}{0}

We investigate a reconstruction method of an unknown cavity inside a heat conductor using boundary measurement data. This problem has broad applications in non-destructive testing, particularly for detecting internal defects in thermal systems; see, for instance, \cite{Ciampa-2018, Clemente-2013-34, Nakamura-2015-31}. Mathematically, it is an inverse boundary value problem for the heat equation, formulated as follows. 

Let $\Omega$ be a bounded domain in $\mathbb{R}^d\,(d = 2,\,3)$ with $C^2$ boundary, and $\nu$ denote the unit outward normal vector to $\partial \Omega$. Physically, $\Omega$ is considered as a heat-conductive medium. Further, let $D \Subset \Omega$ be a cavity with Lipschitz boundary, and assume that $\Omega \setminus \overline{D}$ is connected. For simplicity in describing and understanding arguments, we assume that $\Omega$ and $D$ are simply connected. Also, for any set $B\subset \mathbb{R}^d$, $B_T$ stands for $B \times (0,\,T)$ throughout this paper. 

Now, for a prescribed temperature distribution $g\not\equiv 0$ on $(\partial \Omega)_T$, consider the corresponding temperature field $u$ which is given as the following initial boundary value problem:
\begin{equation}
\label{u-equation}
  \begin{cases}
(\partial_t - \Delta) u =0 & \mathrm{in}\ (\Omega \setminus \overline{D})_T, \\
u = 0 & \mathrm{on}\ (\partial D)_T, \\
u = g & \mathrm{on}\ (\partial \Omega)_T,\\
u = 0 & \mathrm{at}\ t=0.
\end{cases}
\end{equation}
The well-posedness of the forward problem \eqref{u-equation} is well established, as discussed in \cite{Costabel-1990-13, Evans-2010}, for example. 

The inverse boundary value problem of our interest is described as follows.

\medskip
\noindent {\bf Inverse Problem:} Reconstruct the cavity $D$ using a single set of Cauchy data $\{u|_{(\partial\Omega)_T},\,\partial_{\nu} u|_{(\partial \Omega)_T}\}$.

Theoretical results on the uniqueness and conditional stability of this inverse problem are well established (see \cite{Kress-1998-14, Vessella-2008-24} and the references therein). From a numerical perspective, various reconstruction methods such as Newton's methods \cite{Bryan-2005-21, Kress-1998-14}, optimization methods \cite{H-T-2013} and qualitative methods \cite{Nakamura-2009-211, Ikehat-2009-25, Isakov-2010-9, Wang-2018-369, Wang-2012-28, Isozaki-2012-6} are known. Our focus is on a domain sampling method known as the no-response test (NRT). Each sampling domain $G$ called a test domain has to be a simply connected $C^2$-domain (i.e., the boundary $\partial G$ of $G$ is $C^2$-smooth) satisfying $G\Subset\Omega$ (i.e., $\overline G\subset\Omega$) and the connectivity of $\Omega\setminus\overline{G}$. This method has been applied to solve inverse scattering problems \cite{Nakamura-2008-187, Nakamura-2006-31, Potthast-2003-63, Potthast-2007-38, Potthasta-2010-234}, the Oseen problem \cite{Potthast-2016-304} and the inverse boundary value problem for the Laplace equation \cite{Nakamura-2021-15, Sun-2023-485}. In our previous work \cite{Sun-2024-40}, we established the theoretical foundation for using the NRT to solve our inverse boundary value problem through its duality with the range test method (RT), and conducted some numerical experiments. Specifically, for a test domain $G\Subset \Omega$, we designed an indicator function for the range test, denoted as $I_{RT}(G)$. By the analytical extension property of the solution to \eqref{u-equation}, we proved that $I_{RT}(G) < \infty$ if and only if $\overline D\subset \overline G$. This implies the convergence result for the RT. That is, roughly speaking, $D$ is reconstructed by intersecting all test domains containing $D$. We also defined the indicator function $I_{NRT}(G)$ for the NRT and showed that these two methods are dual by proving that $I_{RT}(G) = I_{NRT}(G)$. Thus, we also have the convergence result for the NRT. In the above mentioned papers on the NRT, the convergence results were shown by using the duality relation between the RT and NRT. Hence, as a continuation of our work, we aim to investigate the convergence of the NRT without using the duality.

The main contributions of our study are twofold. First, for a prescribed test domain $G\Subset \Omega$, we define the indicator function $I_{NRT}(G)$ for the NRT and provide an equivalent form. Second, by using the analytic extension property of solutions to the heat equation, we rigorously prove the convergence result which consists of the following statements:
\begin{itemize}
\item [{\rm (i)}] $I_{NRT}(G)<\infty$ if and only if $\overline D \subset \overline G$, provided that $u$ does not admit any analytical extension across $(\partial D)_T$. Consequently, if we call a test domain $G$ positive when $I_{NRT}(G)<\infty$, then this fundamental result enables the reconstruction of $D$ through the intersection of all positive test domains. 
\item[{\rm(ii)}] Additionally, for the case where $\overline D\not\subset\overline G$, we provide a blow-up estimate for $I_{NRT}(G)$ from below.
\end{itemize}

Statement (ii) could be advantageous in improving the numerical algorithms for applying the NRT. We also highlight two advantages of the NRT: (1) it depends on only one boundary measurement; (2) it requires little {\it a priori} information about the cavity.

The rest of this paper is organized as follows. In Section \ref{preliminary}, we introduce some preliminaries. In Section \ref{main-results}, we establish the convergence of the NRT using the analytic extension property. The final section is for conclusions and discussions.

\section{Preliminaries}\label{preliminary}
\setcounter{equation}{0}

In this section, we present some preliminaries, including the notations, definitions and important lemmas used throughout this paper. We begin by introducing the real anisotropic Sobolev space (\cite{Costabel-1990-13}):
$$
H^{r,s}(\mathbb{R}^d \times \mathbb{R}):=L^2(\mathbb{R};H^r(\mathbb{R}^d))\cap H^s(\mathbb{R};L^2(\mathbb{R}^d)), \quad r,\,s \geq 0.
$$
For $r,s <0$, $H^{r,s}(\mathbb{R}^d \times \mathbb{R})$ is defined as the dual space of $H^{-r,-s}(\mathbb{R}^d \times \mathbb{R})$. For any bounded domain $A\subset\mathbb R^d$, we denote by $H^{r,s}\left(A_T\right)$ the space of restrictions of elements in $H^{r,s}(\mathbb{R}^d \times \mathbb{R})$ to $A_T$. The space $H^{r,s}\left((\partial A)_T\right)$ is defined analogously. Moreover, we need the space
$$
\tilde H^{r,s}\left(A_T\right):= \{u(x,t)\in  H^{r,s}\left(A\times (-\infty, T)\right):\,u(x,t) = 0\;\textrm{for}\; t <0\}.
$$ 

Now, we specify the definition of homotopy.
\begin{defn}
A family $\{G_\ell\}_{\ell\in[0,1]}$ of test domains is called a homotopy connecting $\Omega$ and $G$ if it satisfies the following properties:
\begin{itemize}
  \item $G_0 = \Omega$, $G_1 = G$;
  \item Each $G_\ell$ satisfies the same topological property and regularity as those for $G$;
  \item $G_\ell$ depends continuously on $\ell\in[0,1]$;
  \item $G_{\ell_2}\Subset G_{\ell_1}$ for any $0\leq \ell_1<\ell_2\leq 1$.
\end{itemize}
\end{defn}
Denote by $\Gamma(x,t;y,s)$ the Green function of the heat equation in $\Omega_T$ with the homogeneous Dirichlet boundary condition on $(\partial \Omega)_T$. For a test domain $G$, we define a boundary integral operator
\begin{equation*}
 H^{-1/2,-1/4}((\partial\Omega)_T)\ni R[\varphi](x,t) := \int^t_0\int_{\partial G} \partial_{\nu(x)}\Gamma(x,t;y,s)\,\varphi(y,s)\,d\sigma(y)ds,\quad (x,t) \in (\partial\Omega)_T,
\end{equation*}
where $\varphi\in H^{-1/2,-1/4}((\partial G)_T)$ and $\nu(x)$ denotes the unit outward normal vector to $\partial\Omega$ at $x\in\partial\Omega$. Henceforth, we let $X = H^{-1/2,-1/4}((\partial G)_T)$ and $Y = H^{-1/2,-1/4}((\partial \Omega)_T)$ for simplicity. The dual spaces of $X$ and $Y$ are denoted as $X^\prime = H^{1/2,1/4}((\partial G)_T)$ and $Y^\prime = H^{1/2,1/4}((\partial \Omega)_T)$, respectively. Moreover, let $R^*:Y \to X$ and $R^\prime: Y^\prime \to X^\prime$ be the adjoint operator and dual operator of $R$, respectively. It is known that $R^* = \mathcal{I}_X^{-1}R^\prime\mathcal{I}_Y$, where $\mathcal{I}_X : X \to X^\prime$ and $\mathcal{I}_Y : Y \to Y^\prime$ are two isometric isomorphisms satisfying
\begin{align*}
    (\mathcal{I}_X\varphi_1)(\varphi_2) &= (\varphi_1, \varphi_2)_X, \quad \varphi_1, \varphi_2\in X,\\
    (\mathcal{I}_Y\psi_1)(\psi_2) &= (\psi_1, \psi_2)_Y, \quad \psi_1, \psi_2\in Y.
\end{align*}
Here, $(\cdot, \cdot)_X$ and $(\cdot, \cdot)_Y$ are the real inner products in $X$ and $Y$, respectively.
Hence, we have for $\zeta\in Y$ that
\begin{equation}\label{R-dual-adjoint}
    \Vert R^*\zeta\Vert_X = \Vert \mathcal{I}^{-1}_XR^\prime \mathcal{I}_Y \zeta\Vert_X = \Vert R^\prime \tilde \zeta\Vert_{X^\prime}\quad \text{with}\,\,\tilde \zeta = \mathcal{I}_Y \zeta.
\end{equation}

Next, we describe the definition of the analytic extension for solutions to the heat equation and present a relevant lemma.

\begin{defn}\label{AE}
  Let $\Omega$ and $D$ be two open bounded domains as described in Section \ref{Introduction}. We say a solution $u(x,t) \in \tilde H^{1,\frac{1}{2}}((\Omega\setminus \overline{D})_T)$ of the heat equation is analytically extended across $(\partial D)_T$ if and only if there exist a domain $D^{*}$ with Lipschitz boundary such that
  $$
  \overline{D^{*}} \subsetneqq \overline{D},\qquad \mathcal{H}^{n-1}(\partial D^{*}\cap \partial D) > 0,
  $$
and a function $u^{\prime} \in \tilde H^{1,\frac{1}{2}}((\Omega\setminus \overline{D^{*}})_T)$ such that
  \begin{equation*}
    \begin{cases}
      (\partial_t - \Delta) u^{\prime} =0 & \mathrm{in}\ (\Omega \setminus \overline{D^{*}})_T,\\
       u^{\prime} = 0 & \mathrm{at}\ t=0, \\
       u^{\prime} = u & \mathrm{in}\ (\Omega \setminus \overline D)_T,
    \end{cases}
  \end{equation*}
where $\mathcal{H}^{n-1}(\cdot)$ represents the $(n-1)$-dimensional Hausdorff measure.
\end{defn}

\begin{lemma}[\cite{Potthast-2007-38}]\label{Coffecient-bounded}
  Let $Z\Subset \Omega$ be a simply connected $C^2$-domain and $Z_e := \Omega\setminus\overline Z$. Assume that a function $\omega$ satisfies the following two properties:
  \begin{itemize}
    \item $\omega$ is analytic in $Z_e$;

    \item For some small $\rho>0$, the set
    \begin{equation}\label{C-C}
    \left\{\mathcal{C}_q(x) := \sup_{|h|=1}\frac{|(h\cdot\nabla)^q\omega(x)|}{q!}\rho^q,\,q\in \mathbb{Z}_+\right\}\;\text{with ${\mathbb Z}_+:={\mathbb N}\cup\{0\}$}
    \end{equation}
  associated with the Taylor coefficients of $\omega$ is uniformly bounded for all $x\in Z_e\cap\Omega_{\rho}$, where $\Omega_\rho := \{z\in\Omega: \mathrm{dist}(z,\partial \Omega)>\rho\}$. In other words, there exists a constant $C>0$ such that 
  \begin{equation}\label{C-P}
  \mathcal{P}(x):=\sup_{q\in \mathbb{Z}_+}\{\mathcal{C}_{q}(x)\}\le C,\quad x\in Z_e\cap\Omega_{\rho}.
  \end{equation}
  \end{itemize}
Then $\omega$ can be extended into an open neighborhood of $\overline Z_e$. More specifically, there exists a domain $Z^\prime$ with $Z^\prime\Subset Z$ such that $\omega$ admits an analytic extension into the whole $\Omega \setminus \overline{Z^\prime}$, which we simply say $\omega$ is extensible into $\Omega \setminus \overline{Z^\prime}$.
\end{lemma}

\begin{remark}
We highlight that the solution $u(x,t)$ with $(x,t)\in(\Omega\setminus\overline D)_T$ of \eqref{u-equation} is analytic with respect to the spatial variables; see, for instance, \cite{Evans-2010}. If we replace $\omega(x)$ in \eqref{C-C} and \eqref{C-P} with $u(x,t)$, then the second condition of Lemma \ref{Coffecient-bounded} has to be replaced by \begin{equation}\label{C-P'}
\mathcal{P}(x,t):= \sup_{q\in \mathbb{Z}_+}\{\mathcal{C}_{q}(x, t)\}\le C,\quad (x,t)\in (Z_e\cap\Omega_{\rho})_T
\end{equation}
with
\begin{equation*}
\left\{\mathcal{C}_q(x,t) := \sup_{|h|=1}\frac{|(h\cdot\nabla_x)^qu(x,t)|}{q!}\rho^q,\,q\in \mathbb{Z}_+\right\},
\end{equation*}
where $\nabla_x$ denotes the gradient operator for the spatial variable. For simplicity, we still use the same notation $\mathcal{P}(x)$ instead of $\mathcal{P}(x,t)$ in the following analysis.
\end{remark}

Finally, for later use, we introduce an integral operator $K$ given as
\begin{equation}\label{K-operator}
  K[\phi](x,t) = \int^T_t\int_{\partial \Omega} \partial_{\nu(y)}\Phi(y,s;x,t)\,\phi(y,s)\,d\sigma(y)ds, \quad (x,t)\in \Omega_T,
\end{equation}
where 
\begin{equation}\label{Phi}
\Phi(y,s;x,t) := \frac{1}{(4\pi(s-t))^{d/2}}\exp\left(-\frac{|y-x|^2}{4(s-t)}\right)\,\,
\text{with}\,\,\Phi(y,s;x,t) = 0\,\,\text{for}\,\,s\leq t.
\end{equation}

\begin{prop}\label{K-dense}
     Let $E\Subset \Omega$ be a simply connected domain with $C^2$ boundary $\partial E$. Consider $K$ as an operator $K: H^{1/2,1/4}((\partial\Omega)_T)\to H^{1/2,1/4}((\partial E)_T)$. Then $K$ is compact and has a dense range.
 \end{prop}
 {\bf Proof.} The compactness of $K$ is clear. So, we next prove that $K$ has a dense range. The proof is based on the proof given in \cite{Wang-2012-28} for a different setting. Since $K$ is compact, it is enough to prove that the adjoint operator $K^*: H^{1/2,1/4}((\partial E)_T) \to H^{1/2,1/4}((\partial\Omega)_T)$ is injective. Furthermore, it is sufficient to prove $\varphi=0$ if $K^\prime[\varphi]=0$, where $K^\prime: H^{-1/2,-1/4}((\partial E)_T)\to H^{-1/2,-1/4}((\partial\Omega)_T)$ is the dual operator of $K$. By elementary calculations, we have
 \begin{equation*}\label{K*-operator}
   K^\prime[\varphi](x,t) = \int^t_0\int_{\partial E} \partial_{\nu(x)}\Phi(x,t;y,s)\,\varphi(y,s)\,d\sigma(y)ds,\quad (x,t)\in (\partial\Omega)_T,
 \end{equation*}
where $\nu(x)$ is the unit outward normal vector to $\partial \Omega$. For later use, we denote by $\mathcal N(K^\prime)$ the null space of $K^\prime$. 

For $\varphi \in \mathcal N(K^\prime)$, we assume $\varphi$ is such that $\varphi \in C\left((\partial E)_T\right)$ with $\varphi = 0$ at $t=T$. This is possible to assume because the set consisting of such a $\varphi$ is dense in $\mathcal{N}(K^\prime)$ due to $\Phi(x,t;y,s)=0$ for $t\le s$. Then, we continuously extend $\varphi$ to $t>T$ and retain the notation $\varphi$ for simplicity. Define 
\[
v(x,t) = \int^t_0\int_{\partial E} \Phi(x,t;y,s)\,\varphi(y,s)\,d\sigma(y)ds,\quad (x,t)\in (\mathbb R^d \setminus\partial E)\times (0,+\infty).
\]
Note that $v$ satisfies
  \begin{equation*}
    \begin{cases}
      (\partial_t - \Delta) v =0 & \mathrm{in}\ (\mathbb R^d \setminus\overline \Omega)\times (0,+\infty),\\
      \partial_\nu v = 0 & \mathrm{on}\ \partial \Omega\times (0,+\infty),\\
       v = 0 & \mathrm{at}\ t=0.
    \end{cases}
  \end{equation*}

Next, we prove $v = 0$ in $(\mathbb R^d \setminus\overline E)\times (0,+\infty)$. By Lemma 6.16 and the inequality (9.17)
\[
r^\alpha\exp(-r)\leq \alpha^\alpha\exp(-\alpha), \quad 0<r, \alpha<\infty,
\]
both given in \cite{Kress-2014}, we have for $t>s$ and $x\neq y$ that
\begin{equation}\label{E1}
\begin{split}
    |\Phi(x,t;y,s)| &= \frac{1}{(4\pi(t-s))^{d/2}}\exp\left(-\frac{|x-y|^2}{4(t-s)}\right)\\
    &\leq \frac{C_1}{(t-s)^{\beta_1}|x-y|^{d-2\beta_1}},
\end{split}
\end{equation}
\begin{equation}\label{E2}
\begin{split}
    |\partial_{\nu(x)}\Phi(x,t;y,s)| &= \left|\frac{1}{(4\pi(t-s))^{d/2}}\frac{(\nu(x), y-x)}{2(t-s)}\exp\left(-\frac{|x-y|^2}{4(t-s)}\right)\right|\\
    &\leq  \frac{C_2|x-y|^2}{(t-s)^{d/2+1}}\exp\left(-\frac{|x-y|^2}{4(t-s)}\right) \\
    &\leq \frac{C_3}{(t-s)^{\beta_2}|x-y|^{d-2\beta_2}},
\end{split}
\end{equation}
where $0<\beta_1,\,\beta_2<\frac{1}{2}$ and $C_1,\,C_2,\,C_3$ are positive constants. Let $B_R$ be a ball centered at the origin with sufficiently large radius $R$, and denote $\mathcal O_R: = (\mathbb R^d \setminus\overline \Omega)\cap B_R$.
Then we can deduce from \eqref{E1} that $\Vert v(\cdot,+\infty)\Vert_{L^2(\mathcal O_R)} = 0$. Moreover, note that $\Vert v(\cdot,0)\Vert_{L^2(\mathcal O_R)} = 0$ and $\partial_\nu v = 0$ on $\partial \Omega\times (0,+\infty)$. Hence, we have
\begin{equation*}
0 = \int^\infty_0\int_{\mathcal O_R} (v\,\partial_tv - v\,\Delta v)\,d\sigma(x)dt = -\int^\infty_0\int_{\partial B_R} v\,\partial_\nu v\,d\sigma(x)dt + \int^\infty_0\int_{\mathcal O_R} (\nabla v)^2\,d\sigma(x)dt. 
\end{equation*}
Here, note that by \eqref{E1} and \eqref{E2}, we have
\[
\int^\infty_0\int_{\partial B_R} v\,\partial_\nu v\,d\sigma(x)dt\to 0\,\, \mathrm{as}\,\,R\to +\infty.
\]
Hence, we have
\[
\lim_{R\to +\infty}\int^\infty_0\int_{\mathcal O_R} (\nabla v)^2\,d\sigma(x)dt = 0.
\]
Taking account of $v \to 0$ on $(\partial B_R)\times(0,+\infty)$ as $R\to +\infty$,
this implies $v = 0$ in $\mathcal O_R\times (0,+\infty)$. Then, by the unique continuation property of solutions to the heat equation, we have $v = 0$ in $(\mathbb R^d \setminus\overline E)\times (0,+\infty)$.

For brevity, we use `+' and `-' to indicate the limits taken from the exterior and interior of $E$. The above argument allows us to derive that $v^+|_{(\partial E)_T} = 0$ and $\partial_\nu v^+|_{(\partial E)_T} = 0$, which implies $v^-|_{(\partial E)_T} = 0$ by the jump relation of the single-layer potential. Hence, $v^-$ satisfies
  \begin{equation*}
    \begin{cases}
      (\partial_t - \Delta) v^- =0 & \mathrm{in}\ E_T,\\
      v^- = 0 & \mathrm{on}\ (\partial E)_T,\\
      v^- = 0 & \mathrm{at}\ t=0.
    \end{cases}
  \end{equation*}
It follows that $v=0$ in $E_T$ by the uniqueness of solutions to this problem, thereby $\partial_\nu v^-|_{(\partial E)_T} = 0$. Using the jump relations of the layer potentials for the heat equation, we have $\varphi = 0$ by observing $\partial_\nu v^+|_{(\partial E)_T} = \partial_\nu v^-|_{(\partial E)_T} = 0$. The proof is complete.
$\hfill\qedsymbol$

\section{Main results}\label{main-results}
\setcounter{equation}{0}

In this section, we show the convergence of the NRT by means of the analytic extension property of solutions to the heat equation.

For a given $g\in Y^\prime$, let $\mu \in \tilde H^{1,\frac{1}{2}}(\Omega _T)$ be the background solution that satisfies
\begin{equation}\label{mu-equation}
  \begin{cases}
(\partial_t - \Delta) \mu =0 & \mathrm{in}\ \Omega_T, \\
\mu = g & \mathrm{on}\ (\partial \Omega)_T,\\
\mu = 0 & \mathrm{at}\ t=0.
\end{cases}
\end{equation}
Assume that $u\in \tilde H^{1,\frac{1}{2}}(\Omega \setminus \overline{D})_T$ is the solution to \eqref{u-equation}. Define $w:= u-\mu$. Then we have
\begin{equation}
\begin{cases}\label{Omega-equation}
(\partial_t - \Delta) w =0 & \mathrm{in}\ (\Omega \setminus \overline{D})_T, \\
w = -\mu & \mathrm{on}\ (\partial D)_T, \\
w = 0 & \mathrm{on}\ (\partial \Omega)_T,\\
w = 0 & \mathrm{at}\ t=0.
\end{cases}
\end{equation}

We now consider the following boundary integral equation:
\begin{equation}\label{main-equation}
  R[\varphi](x,t) = \partial_{\nu} w(x,t),\quad (x,t) \in (\partial\Omega)_T.
\end{equation}

\begin{lemma}[\cite{Sun-2024-40}]\label{Analytic}
  The equation \eqref{main-equation} is solvable in $X$ if and only if $v$ can be analytically extended into $(\Omega\setminus\overline G)_T$ and $w_+|_{(\partial G)_T} \in  X^\prime$, where $w_+ = w |_{(\Omega\setminus \overline{G})_T}$.
\end{lemma}

For a test domain $G$, we recall the indicator function of the NRT (see \cite{Sun-2024-40}):
\begin{equation}\label{ind-NRT}
  I_{NRT}(G) = \sup_{\zeta \in Y, \,\Vert R^*\zeta\Vert_X \leq 1} |(\zeta,\partial_\nu w)_Y|.
\end{equation}
We demonstrate that $I_{NRT}(G)$ has an equivalent form in the following proposition.
\begin{prop}
  The indicator function $I_{NRT}(G)$ defined by \eqref{ind-NRT} is equivalent to
  \begin{equation}\label{NRT-F2}
    I_{NRT}(G) = \sup_{\phi \in Y^\prime, \,\Vert R^\prime K[\phi]\Vert_{X^\prime} \leq 1} \left|\int^T_0\int_{\partial D}\left( K[\phi](y,s)\,\partial_\nu w(y,s) - \partial_\nu K[\phi](y,s)\,w(y,s)\right)\, d\sigma(y)ds \right|.
  \end{equation}
\end{prop}

{\bf Proof.} For $\zeta\in Y$, we have from the definition of $\mathcal{I}_Y$ that
\begin{equation*}
  (\zeta, \partial_\nu w)_Y = (\mathcal{I}_Y\zeta)(\partial_\nu w) = \tilde\zeta(\partial_\nu w) = \int^T_0\int_{\partial \Omega}\tilde{\zeta}\,\partial_\nu w \,d\sigma(y)\,ds,
\end{equation*}
where the last integral is in the dual sense. The combination of \eqref{R-dual-adjoint} and \eqref{ind-NRT} yields that
\[
I_{NRT}(G) = \sup_{\tilde \zeta \in Y^\prime, \,\Vert R^\prime\tilde\zeta\Vert_{X^\prime} \leq 1} \left|\int^T_0\int_{\partial \Omega}\tilde{\zeta}\,\partial_\nu w \,d\sigma(y)\,ds\right|.
\]
Note that $K[\phi]$ admits the following jump relation (see, for instance, \cite{Wang-2012-28}):
\begin{eqnarray*}
  K[\phi](x,t) &=& \int^T_t\int_{\partial \Omega} \partial_{\nu(y)}\Phi(y,s;x,t)\,\phi(y,s)\,d\sigma(y)ds - \frac{1}{2}\phi(x,t)\\
  &=:& \mathcal K[\phi](x,t)- \frac{1}{2}\phi(x,t), \quad (x,t)\in (\partial \Omega)_T,
\end{eqnarray*}
where $\mathcal K: Y^\prime\to Y^\prime$ is a compact operator.
Furthermore, we observe that $K[\phi](x,t)$ satisfies the homogeneous backward heat equation in $\Omega_T$ and $K[\phi](x,T) = 0$ for $x\in \Omega$.
Combining the uniqueness of solutions to the backward heat equation and the Fredholm alternative theorem, we know that for any $\tilde{\zeta} \in Y^\prime$, there exists a unique $\phi \in Y^\prime$ such that
\begin{equation*}
  K[\phi](x,t) = \tilde{\zeta}(x,t),\quad (x,t)\in (\partial \Omega)_T.
\end{equation*}
As a result, we are led to
\begin{equation*}
    I_{NRT}(G) = \sup_{\phi \in Y^\prime, \,\Vert \mathcal R^\prime K[\phi]\Vert_{X^\prime} \leq 1} \left|\int^T_0\int_{\partial \Omega} K[\phi](y,s)\,\partial_\nu w(y,s)\, d\sigma(y)ds \right|.
\end{equation*}

Integrating by parts over $(\Omega\setminus \overline D)_T$, we deduce that
\begin{equation*}
\begin{split}
   & \int^T_0\int_{\partial \Omega} K[\phi](y,s)\,\partial_\nu w(y,s) \,d\sigma(y)ds \\
    = & \int^T_0\int_{\partial \Omega}\left( K[\phi](y,s)\,\partial_\nu w(y,s) - \partial_\nu K[\phi](y,s)\,w(y,s)\right)\, d\sigma(y)ds\\
    = & \int^T_0\int_{\partial D}\left( K[\phi](y,s)\,\partial_\nu w(y,s) - \partial_\nu K[\phi](y,s)\,w(y,s)\right)\, d\sigma(y)ds,
\end{split}
\end{equation*}
and therefore \eqref{NRT-F2} holds. The proof is complete. $\hfill\qedsymbol$
 
Now, we are in a position to state our main result, which reveals the convergence of the NRT.

\begin{theorem}\label{Main-theorem}
Let $\Omega$ be a bounded domain in $\mathbb{R}^d\,(d = 2,\,3)$ with $C^2$ boundary, and $D \Subset \Omega$ be a cavity with Lipschitz boundary. Assume that $u$ is the solution to \eqref{u-equation}. For a test domain $G$, we have the following results:
\begin{itemize}
  \item[(1)] If $\overline{D}\subset \overline{G}$, then $I_{NRT}  (G) <\infty$.

  \item[(2)] If $\overline{D}\not\subset \overline{G}$ but $u$ is extensible into $(\overline D\setminus G)_T$ with Lipschitz smooth boundary $\partial(\overline D\setminus G)$ of $\overline D\setminus G$, then $I_{NRT}(G) <\infty$.

  \item[(3)] If $\overline {D}\not\subset \overline{G}$ and $u$ cannot be analytically extended across $(\partial D)_T$, then $I_{NRT}(G) =\infty$.
\end{itemize}
\end{theorem}

{\bf Proof.} For the simplicity of writing the figures in the proof, we restrict our argument to the case of $d=2$, but the proof works without any change for $d=3$. Let us start with the case (1), i.e., $\overline{D}\subset \overline{G}$. We know from Definition \ref{AE} that $u$ admits an analytic extension into $G_T$, and therefore $u|_{(\partial G)_T} \in X^\prime$, which implies $w|_{(\partial G)_T} \in X^\prime$ as well. According to Lemma \ref{Analytic}, there exists a solution $\varphi\in X$ satisfying \eqref{main-equation}. Hence,
\begin{equation*}
  I_{NRT}(G) = \sup_{\zeta\in Y,\,\Vert R^*\zeta\Vert_X\leq 1}|(\zeta, R\varphi)_{Y}|= \sup_{\zeta\in Y,\,\Vert  R^*\zeta\Vert_X\leq 1}|(R^*\zeta,\varphi)_{X}|< \infty.
\end{equation*}

For the case (2), we also have $w|_{(\partial G)_T} \in X^\prime$ by the assumption. Then, we can derive $I_{NRT}(G) < \infty$ using the same argument as in the case (1).

\begin{figure}[htbp]
 \centering
  \subfigure[]
  {
	\begin{minipage}{6cm}
	\centering
	\includegraphics[scale=0.28]{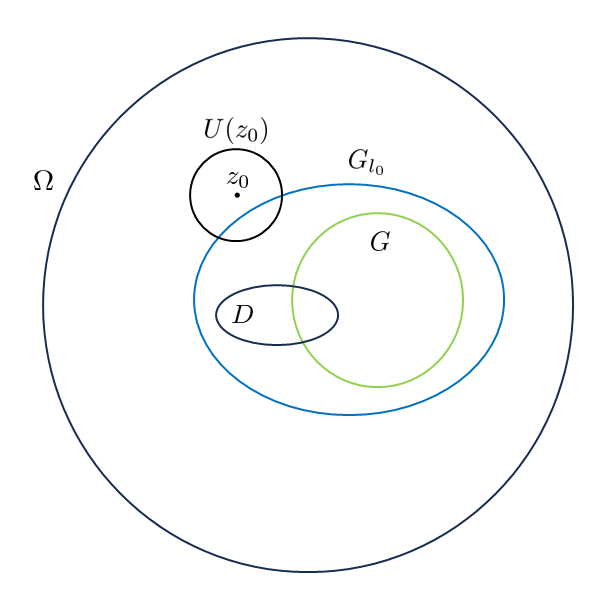}
	\end{minipage}
  }
  \subfigure[]
  {
	\begin{minipage}{6cm}
	\centering
	\includegraphics[scale=0.35]{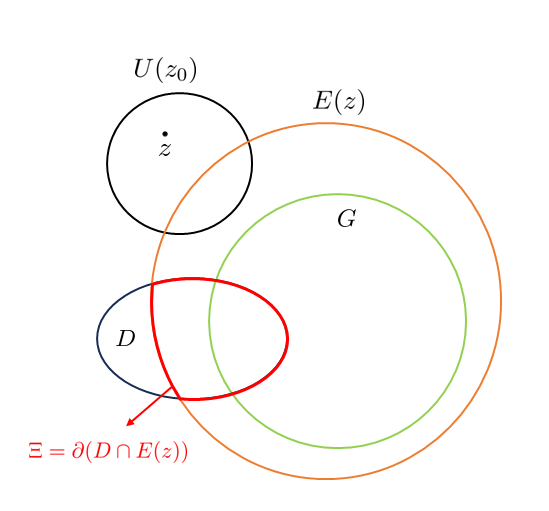}
	\end{minipage}
  }
  \caption{(a) The positional relationship among $G_{l_0}$, $U(z_0)$, $D$ and $G$. (b) The positional relationship among $U(z_0)$, $E(z)$, $D$ and $G$. }
 \label{illustration}
\end{figure}

Next, we focus on the case (3), for which the proof is more complicated than the previous two cases. Let $\{G_l\}_{l\in[0,1]}$ be a homotopy connecting $\Omega$ and $G$.
Since $u$ cannot be analytically extended across $(\partial D)_T$ and $\overline {D}\not\subset \overline{G}$, there exists $l_0\in(0,1)$ with the property that $u$ is extensible into $(\Omega \setminus G_l)_T$ for any $l<l_0$ but not for all $l> l_0$. By the definition of the homotopy, we have $G\Subset G_{l_0}$. Also, according to Lemma \ref{Coffecient-bounded}, $\mathcal{P}(z)$ cannot be uniformly bounded for $z\in (\Omega\setminus \overline{G_{l_0}})\cap\Omega_{\rho}$ for some $\rho>0$. That is, there exists at least one $z_0\in \Omega\setminus \overline {G_{l_0}}$ with dist$(z_0,\partial\Omega)>\rho$ such that $\mathcal{P}(z)$ cannot be uniformly bounded for any neighborhood $U(z_0)$ of $z_0$. See Figure \ref{illustration} (a) for the positional relationship among $G_{l_0}$, $D$ and $G$ which we have here.

For the above $z_0$ and $\rho_G := \text{dist}(z_0,\overline{G})$, take an open neighborhood $U(z_0)\subset B(z_0,\rho_G/2)$ of $z_0$, where $B(z_0,\rho_G/2)$ is a ball centered at $z_0$ with radius $\rho_G/2$. Fix any $h\in{\mathbb R}^2$ with $|h|=1$, define 
\begin{equation*}
  N_{(m;h)}(z,s):=\Vert (h\cdot \nabla_z)^m\Phi(z,s; \cdot,\cdot)\Vert_{H^{1,1/2}\left(\left(G\cup\left(\Omega\setminus\overline\Omega_\varepsilon\right)\right)_T\right)}
\end{equation*}
for $m\in \mathbb{Z}_+,\,z\in W(z_0):= U(z_0)\cap (\Omega_\varepsilon\setminus \overline G_{l_0}),\,s\in(0,T),$ where $\Omega_\varepsilon:= \{z: \text{dist}(z,\partial \Omega)>\varepsilon\}$ with some small constant $\varepsilon>0$. Let
\begin{equation*}
  \mathcal{V}^{(z,s)}_{(m,h)}(x,t):=\frac{1}{2\Vert \mathcal{T} \Vert\Vert R^\prime \Vert N_{(m;h)}(z,s)}(h\cdot \nabla_z)^m\Phi(z,s;x,t),
\end{equation*}
where $\mathcal{T}: \tilde H^{1,1/2}(G_T)\to H^{1/2,1/4}((\partial G)_T)$ is the trace operator. Then, we clearly have
\[
\Vert\mathcal{V}^{(z,s)}_{(m,h)}\Vert_{H^{1,1/2}(G_T)}\leq \frac{1}{2\Vert\mathcal{T} \Vert \Vert \mathcal R^\prime \Vert}.
\]

Moreover, for any fixed $z\in W(z_0)$, there is a simply connected domain $E(z)\subset \Omega$ satisfying the following properties (see Figure \ref{illustration} (b)):
\begin{itemize}
  \item[(1)] $z\not\in E(z)\,\text{and}\,\mathrm{dist}(z, E(z)) > 0$;
  \item[(2)] $G\subset E(z)$;
  \item[(3)] The boundary $\partial E(z)$ of $E(z)$ is of class $C^2$;
  \item[(4)] $u$ satisfies the homogeneous heat equation in $(\Omega\setminus E(z))_T$.
\end{itemize}
Then we can easily confirm the following properties:
\begin{itemize}
\item[{\rm(i)}] $E(z)\cap D \neq \emptyset$; 
\item[{\rm(ii)}] $\mathcal{V}^{(z,s)}_{(m,h)}(x,t)\in H^{1/2,1/4}((\partial E(z))_T)$; 
\item[{\rm(iii)}] $K: H^{1/2,1/4}((\partial\Omega)_T)\to H^{1/2,1/4}((\partial E(z))_T)$ has a dense range (see Proposition \ref{K-dense}).
\end{itemize}
From the property (iii), there exists a sequence $\{\phi^{z,s,m,h}_n \in H^{1/2,1/4}((\partial\Omega)_T),\,n=1,2,\cdots\}$ such that
\begin{equation*}
  \Vert K[\phi^{z,s,m,h}_n] - \mathcal{V}^{(z,s)}_{(m;h)}\Vert_{H^{1/2,1/4}((\partial E(z))_T)} \to 0,\,\,n\to \infty.
\end{equation*}
Since both $K[\phi^{z,s,m,h}_n](x,t)$ and $\mathcal{V}^{(z,s)}_{(m;h)}(x,t)$ satisfy the homogeneous backward heat equation in $(E(z))_T$ and $K[\phi^{z,s,m,h}_n](\cdot,T) = \mathcal{V}^{(z,s)}_{(m;h)}(\cdot,T) = 0$ in $E(z)$, we have
\begin{equation}\label{con-K}
\Vert K[\phi^{z,s,m,h}_n] - \mathcal{V}^{(z,s)}_{(m;h)}\Vert_{H^{1,1/2}((E(z))_T)} \to 0,\,\, n\to \infty.
\end{equation}
Therefore, we can find a subsequence $\{\phi^{z,s,m,h}_{n_k} \in H^{1/2,1/4}((\partial\Omega)_T),\,k=1,2,\cdots\}$ such that
\begin{equation*}
\Vert K[\phi^{z,s,m,h}_{n_k}]-  \mathcal{V}^{(z,s)}_{(m;h)}\Vert_{H^{1,1/2}(G_T)}\leq\frac{1}{2\Vert \mathcal{T}\Vert \Vert R^\prime \Vert}.
\end{equation*}
Henceforth, we denote $\xi_k:= \phi^{z,s,m,h}_{n_k}$ for simplicity. Then we have
\begin{eqnarray*}
\Vert K[\xi_k]\Vert_{H^{1,1/2}(G_T)}&\leq& \Vert K[\xi_k] - \mathcal{V}^{(z,s)}_{(m;h)}\Vert_{H^{1,1/2}(G_T)}+\Vert \mathcal{V}^{(z,s)}_{(m;h)}\Vert_{H^{1,1/2}(G_T)} \\
  &\leq& \frac{1}{\Vert\mathcal{T}\Vert \Vert R^\prime \Vert},
\end{eqnarray*}
which implies $\Vert R^\prime K[\xi_k]\Vert_{X^\prime} \leq 1$.

Next, we define
\begin{equation*}
    I^{(k)}_{NRT}(G):= \left|\int^T_0\int_{\partial D} \left(K[\xi_k](x,t)\,\partial_\nu w(x,t) - \partial_\nu K[\xi_k](x,t)\,w(x,t)\right)\,d\sigma(x)dt \right|.
\end{equation*}
Recall that $w$ and $K[\xi_k]$ satisfy the homogeneous forward and backward heat equation in $(\Omega\setminus \overline{E(z)})_T$, respectively. Applying Green's theorem over $(D\setminus \overline{E(z)})_T$, we have
\begin{equation*}
    I^{(k)}_{NRT}(G) = \left|\int^T_0\int_{\Xi} \left(K[\xi_k](x,t)\,\partial_\nu w(x,t) - \partial_\nu K[\xi_k](x,t)\,w(x,t)\right)\,d\sigma(x)dt \right|,
\end{equation*}
where $\Xi:= \partial (D\cap E(z))\subset \overline {E(z)}$. Then we deduce from \eqref{con-K} that 
\begin{eqnarray*}
  \lim_{k\to \infty}I^{(k)}_{NRT}(G) &=& \lim_{k\to \infty} \left|\int^T_0\int_{\Xi} \left(K[\xi_k](x,t)\,\partial_\nu w(x,t) - \partial_\nu K[\xi_k](x,t)\,w(x,t)\right)\,d\sigma(x)dt \right|\\
  &=& \left|\int^T_0\int_{\Xi}\left(\mathcal{V}^{(z,s)}_{(m;h)}(x,t)\,\partial_\nu w(x,t) - \partial_\nu \mathcal{V}^{(z,s)}_{(m;h)}(x,t)\,w(x,t)\right)\,d\sigma(x)dt \right|\\
  &=& \frac{1}{2\Vert\mathcal{T} \Vert \Vert\mathcal R^\prime \Vert N} \Big|\int^T_0\int_{\Xi}\left[(h\cdot \nabla_z)^m\Phi(z,s;x,t)\,\partial_{\nu} w(x,t)\right.  \\
  & &\quad\quad\quad\quad\quad\quad\quad\quad\quad\quad \left. - \partial_{\nu} ((h\cdot \nabla_z)^m\Phi(z,s;x,t))\, w(x,t)\right]\,d\sigma(x)dt \Big|
\end{eqnarray*}
with $N := N_{(m;h)}(z,s)$. By Green's formula for the heat equation, we have
\begin{equation}
\begin{split}
\lim_{k\to \infty}I^{(k)}_{NRT}(G) &= \left|\frac{1}{2\Vert\mathcal{T}\Vert \Vert R^\prime \Vert N} (h\cdot \nabla_z)^m w(z,s)-\frac{1}{2\Vert\mathcal{T} \Vert \Vert \mathcal R^\prime \Vert N} (h\cdot \nabla_z)^m \tilde w(z,s)\right|\\
  &=: \left| I_{(m,h)}(z,s) - \tilde I_{(m,h)}(z,s)\right|,
\end{split}
\end{equation}
where
\begin{equation*}
  \tilde w(z,s):= \int^T_0\int_{\partial \Omega} \left(\Phi(z,s;x,t)\,\partial_{\nu} w(x,t) - \partial_{\nu(x)}\Phi(z,s;x,t)\,w(x,t)\right)\,d\sigma(x)dt.
\end{equation*}
Note that $\mathcal{P}(z)$ is not uniformly bounded for $z\in W(z_0)$. Then, for any $\kappa>0$, we can find $z_\kappa\in W(z_0),\,m_\kappa\in \mathbb{Z}_+, s_\kappa \in (0,T)$ and $|h_\kappa|=1$ such that
\begin{equation}\label{Unbounded-estimate}
\left|\frac{(h_\kappa\cdot \nabla)^{m_\kappa} w(z_\kappa,s_\kappa)}{m_\kappa!} \rho^{m_\kappa}\right|\ge\kappa,\quad \left|\tilde I_{(m_\kappa,h_\kappa)}(z_\kappa,s_\kappa)\right|\le \tilde{c},
\end{equation}
with some positive constant $\tilde{c}$ independent of $\kappa$. Furthermore, there exists a positive constant $c$ which satisfies
\begin{equation}\label{Inequality-Gamma}
2\Vert\mathcal{T} \Vert \Vert\mathcal R^\prime \Vert N_{(m;h)}(z,s) \le c\frac{m!}{\rho^m},\,\,\,z\in W(z_0),\,s\in(0,T),\,m\in \mathbb{Z}_+,\, |h|=1.
\end{equation}
Combining \eqref{Unbounded-estimate} and \eqref{Inequality-Gamma}, we have
\begin{equation}
\begin{split}
\left| I_{(m_\kappa,h_\kappa)}(z_\kappa,s_\kappa)-\tilde I_{(m_\kappa,h_\kappa)}(z_\kappa,s_\kappa)\right| &\geq \left| I_{m_\kappa,h_\kappa}(y_\kappa,s_\kappa)\right|- \left|\tilde I_{m_\kappa,h_\kappa}(y_\kappa,s_\kappa)\right| \\
&\geq \frac{\kappa}{c}-\tilde{c}\to \infty\,\, \mathrm{as}\,\,\kappa\to \infty.
\end{split}
\end{equation}
Consequently, we have
\begin{equation}
\begin{split}
  I_{NRT}(G) &=\sup_{\phi \in Y^\prime, \,\Vert \mathcal R^\prime K[\phi]\Vert_{X^\prime} \leq 1} \left|\int^T_0\int_{\partial D} \left(K[\phi]\,\partial_\nu w - \partial_\nu K[\phi]\,w\right)\, d\sigma(y)ds \right|\\
  &\geq \lim_{k\to \infty} I^{(k)}_{NRT}(G)\\
  &= \infty.
\end{split}
\end{equation}
Hence, we have completed the proof.
$\hfill\qedsymbol$

\section{Conclusion and discussion}\label{Conclusion}

In this paper, we consider an inverse problem of reconstructing an unknown cavity compactly embedded in a heat-conductive medium by one measurement at the boundary of the medium. Specifically, we showed the convergence of the reconstruction method called the NRT without using the convergence of another reconstruction method via the duality between these two methods. More precisely, using the analytical extension property of solutions to the heat equation, we proved that the indicator function of the NRT for a test domain converges if and only if the test domain fully contains the cavity. This result provides a rigorous criterion for reconstructing the cavity from only a single boundary measurement. 

The proof of Theorem \ref{Main-theorem} shows that the convergence of the NRT completely relies on the analytic continuation of the solution to the forward problem. Even though it is very hard to give a general criterion for the boundary of the cavity creating singularities in the solution $u$ of the forward problem, the situation is different for the numerical realization of the NRT. More precisely, when we discretize the medium with cavity, the boundary of the cavity creates singularities in the solution $u$ of the forward problem. This could be a numerical advantage of the NRT. In fact, our numerical study of the NRT given in \cite{Sun-2024-40} indicates such an observation.

\bigskip
{\bf Acknowledgement:}
The first author was supported by the National Natural Science Foundation of China (Nos. 12471395, 12071072). The second author was supported by the JSPS KAKENHI (Grant Number 22K03366, 25K07076). The third author was supported by the National Natural Science Foundation of China (Nos. 12471395, 12071072, 12241102) and the Jiangsu Provincial Scientific Research Center of Applied Mathematics under Grant No. BK20233002.

\end{document}